\documentclass[a4paper]{article}

\usepackage{delarray,verbatim,enumerate,a4}
\usepackage{amsmath,amsthm,amstext,amsbsy,amssymb,amsfonts,amscd}
\usepackage[all]{xy}
\usepackage[frenchb]{babel}
\usepackage{aeguill}

\newtheorem{Th}{Th\'eor\`eme}[]

\newtheorem{Prop}[Th]{Proposition}
\newtheorem{Cor}[Th]{Corollaire}

\newtheorem{Sco}[Th]{Scolie}
\newtheorem{Def} [Th]{D\'efinition}

\newtheorem{DProp}[Th]{D\'efinition \& Proposition}

\def\Preuve{\smallskip\noindent {\it Preuve.~}}
\def\PreuveTh{\smallskip\noindent {\it Preuve du Théorème.~}}
\def\Remarque{\smallskip\noindent {\it Remarque.~}}

\def\Exemple{\smallskip\noindent {\it Exemple.~}}

\font\teneufm=eufm10
\font\seveneufm=eufm7
\font\fiveeufm=eufm5
\newfam\gothfam
\textfont\gothfam=\teneufm \scriptfont\gothfam=\seveneufm
\scriptscriptfont\gothfam=\fiveeufm
\def\goth{\fam\gothfam}

\def\RR{\mathbb R}				\def\CC{\mathbb C}			\def\QQ{\mathbb Q}	
\def\NN{\mathbb N}			\def\ZZ{\mathbb Z}
\def\F2{\mathbb{F}_2}		\def\Z2{\mathbb{Z}_2}	
\def\Zl{\mathbb{Z}_\ell} 	\def\Ql{\mathbb{Q}_\ell}

 				\def\U{\mathcal  U}
\def\J{\mathcal  J}  				\def\R{\mathcal  R}
 	  	\def\Cl{\mathcal  C \ell}
\def\E{\mathcal  E}		\def\V{\mathcal  V}

		\def\p{{\goth p}}		\def\l{{\goth l}}	\def\L{{\goth L}}

\def\wi{\widetilde}	\def\ov{\overline}	\def\wh{\widehat}
\def\rg{\operatorname{rg}}	\def\Ind{\operatorname{Ind}}
	\def\deg{\operatorname{deg}}
\def\Gal{\operatorname{Gal}}	\def\Log{\operatorname{Log}}

\begin{document}

\title{\bf\Large Plongements $\ell$-adiques et $\ell$-nombres de Weil}

\author{ Jean-François {\sc Jaulent} }
\date{}
\maketitle

{\small
\noindent{\bf Résumé.} Nous introduisons la notion de nombre de Weil $\ell$-adique par analogie avec la notion classique de nombre de Weil à l'infini~; et nous en étudions quelques propriétés en liaison avec les plongements et les valeurs absolues réelles ou $\ell$-adiques des corps de nombres. En appendice, nous en tirons diverses applications à la théorie d'Iwasawa des tours cyclotomiques.
}

\

{\small
\noindent{\bf Abstract.} We define $\ell$-adic analogs of classical Weil numbers in connexion both with complex or $\ell$-adic imbeddings of number fields and real or $\ell$-adic absolute values. As an application we give some consequences related to the Iwasawa theory of cyclotomic towers. 
}

\


\bigskip\smallskip

\noindent{\large\bf Introduction}\medskip

Le but de l'article est de définir les analogues $\ell$-adiques des nombres de Weil et d'étudier quelques unes de leurs propriétés algébriques en liaison avec la théorie d'Iwasawa des corps cyclotomiques.\smallskip

Les nombres de Weil classiques (i.e. à l'infini) interviennent naturellement tant en théorie des nombres qu'en géométrie algébrique ou même en théorie des groupes~: sommes de Gauss ou de Jacobi, tables de caractères, équations fonctionnelles des fonctions $L$ d'Artin, factorisation des fonctions $\zeta$ (de Weil) attachées aux variétés projectives non singulières sur les corps finis, etc (cf. e.g. \cite{BF, GO, La}). C'est d'ailleurs de ce dernier exemple que provient leur nom, en raison de la célèbre conjecture formulée par A. Weil (cf. \cite{We}) et prouvée par R. Deligne (cf. \cite{De}).\smallskip

Du point de vue de l'arithmétique, les nombres de Weil (après normalisation) ne sont rien d'autre que les nombres algébriques dont tous les conjugués réels ou complexes sont de module 1, autrement dit les éléments de la clôture algébrique $\ov\QQ$ de $\QQ$, dont tous les plongements à l'infini sont de module 1.\smallskip

Il nous a donc semblé intéressant de regarder ce que deviennent les propriétés purement algébriques de ces nombres lorsqu'on remplace la place à l'infini par une place finie $\ell$ , à l'instar de ce qui est fait dans \cite{Dn}, ce qui amène ici à considérer $\ov\QQ$ non plus comme un sous-corps de $\CC$, mais comme un sous-corps de $\CC_\ell$ (en fait comme un sous-corps de $\ov\QQ_\ell$) et à remplacer les {\it valeurs absolues réelles} (i.e. à valeurs réelles) par leurs analogues $\ell$-adiques~: {\it les valeurs absolues $\ell$-adiques}.\smallskip

A l'exercice, on constate que si certaines analogies formelles sont bien conservées (par exemple, la notion de corps totalement $\ell$-adique correspond assez naturellement à celle de corps totalement réel), d'autres exigent des restrictions sévères de degré et de ramification (ainsi la notion de conjugaison $\ell$-adique quadratique non ramifiée) et surtout, contrairement à ce qui advient dans le cas classique où les nombres de Weil qui sont des unités (à l'infini) sont des racines de l'unité, l'ensemble des $\ell$-nombres de Weil qui sont des unités (autrement dit le noyau de toutes les valeurs absolues $\ell$-adiques) est généralement un sous-groupe non trivial du groupe multiplicatif du corps considéré, sous-groupe dont le $\ell$-adifié joue un rôle mystérieux dans certaines questions de théorie d'Iwasawa mettant en jeu les conjectures de Gross ou de Leopoldt.

\bigskip \smallskip 
\noindent {\large \bf 1. Corps à conjugaison complexe et nombres de Weil à l'infini} 
\medskip 

Nous revenons dans cette première section sur quelques propriétés galoisiennes des corps à conjugaison complexe et des nombres de Weil classiques (i.e. à l'infini). \par

Quoique la plupart d'entre elles soient essentiellement bien connues (cf. e.g. \cite{GO}), il nous a paru utile, pour la commodité du lecteur, d'en rassembler succintement les preuves. Celles-ci nous servent, en effet de fil directeur dans la section suivante pour définir et étudier les analogues $\ell$-adiques de ces notions.

\medskip
\noindent{\bf 1.a. Corps totalement réels et corps à conjugaison complexe}\medskip

Nous supposons fixés dans cette section un plongement de la clôture algébrique $\ov\QQ$ de $\QQ$ dans le corps des complexes $\CC$~; nous regardons ainsi $\ov\QQ$ comme un sous-corps de $\CC$. Cela permet de parler sans ambiguïté de {\it la} conjugaison complexe $\tau$.

\begin{Def}\label{CM}
 On dit qu'un corps de nombres $K$ est un {\em corps à conjugaison complexe} (ou encore qu'il est de type CM) lorsque c'est une extension quadratique totalement imaginaire d'un sous-corps totalement réel $K_+$.
\end{Def}

En d'autres termes :

\begin{Prop} Les corps de nombres à conjugaison complexe sont ceux de la forme $K = K_+[\sqrt{-\delta}]$, où $K_+$ est un corps totalement réel et $\delta \gg 0$ un élément de $K_+$ totalement positif.
\end{Prop}

Le résultat suivant permet alors de ramener l'étude des corps à conjugaison complexe au cas absolument galoisien~:

\begin{Prop} Le compositum $K = \prod_{i=0}^n K_i$ d'une famille de corps de nombres à conjugaison complexe est encore un corps à conjugaison complexe. \par
En particulier la clôture galoisienne d'un corps à conjugaison complexe est un corps à conjugaison complexe qui est extension quadratique d'un sous-corps galoisien (totalement) réel.
\end{Prop}

\Preuve Pour chaque $i \in \{ 0,\dots,n \}$, écrivons $K_i = F_i[\sqrt{-\delta_i}]$, avec $F_i$ totalement réel et $\delta_i \gg 0$ totalement positif dans $F_i$~; et notons $F = \prod F_i$ le compositum des $F_i$. 
Nous avons immédiatement~: \smallskip

\centerline{$K = F[\sqrt{-\delta_0},\dots,\sqrt{-\delta_n}] =
 F[\sqrt{\delta_1/\delta_0},\dots,\sqrt{\delta_n/\delta_0}] [\sqrt{-\delta_0}]$}\smallskip

\noindent avec $F[\sqrt{\delta_1/\delta_0},\dots,\sqrt{\delta_n/\delta_0}]$ totalement réel et $\delta_0 \gg 0$~; d'où le premier résultat.\smallskip

Pour établir le second , partons de $K_0= F_0[\sqrt{-\delta_0}]$, avec $F_0$ totalement réel et $\delta_0 \gg 0$ totalement positif dans $F_0$, et considérons la famille $K_0^{\sigma_i} = F_0^{\sigma_i}[\sqrt{-\delta_0^{\sigma_i}}]$ des conjugués de $K_0$. Nous obtenons comme plus haut~:
\smallskip

\centerline{$K =  F[\sqrt{\delta_0^{\sigma_1}/\delta_0}, \dots ,
\sqrt{\delta_0^{\sigma_n}/\delta_0}][\sqrt{-\delta_0}]$,}\smallskip

\noindent avec $F = \prod_{i=1}^n F_0^{\sigma_i}$ et $F[\sqrt{\delta
_0^{\sigma_1}/\delta_0}, \dots ,\sqrt{\delta_0^{\sigma_n}/\delta_0}]$ 
galoisien (totalement) réel, engendré sur $F$ par les racines carrées
des $\delta_0^\sigma/\delta_0^{\sigma'}$ pour $\sigma$ et $\sigma'$ 
dans $\Gal(\overline{\mathbb Q}/\mathbb Q)$.\medskip

En particulier, il vient ainsi~:

\begin{Th} Soit $K$ un corps de nombres à conjugaison complexe qui est 
galoisien sur un sous-corps totalement réel $F$ et soit $G$ 
le groupe de Galois $\Gal(K/F)$. Dans ces conditions~:
\begin{itemize} 
\item[(i)] La conjugaison complexe $\tau$ est centrale dans $G$ et 
son corps des points fixes $K_+ = K^{<\tau>}$ est le plus grand sous-corps 
totalement réel de $K$.
\item[(ii)] Pour toute extension quadratique totalement imaginaire
$k = F[\sqrt{-d}]$ de $F$ qui n'est pas contenue dans $K$, l'extension 
composée $K' = Kk = K[\sqrt{-d}]$ est le produit direct de $k$ et d'une
sous-extension galoisienne totalement réelle $K'_+$ de $F$ de 
groupe de Galois isomorphe à $G$.
\end{itemize}
\end{Th}

\Remarque Il suit de là que toute extension absolument galoisienne
à conjugaison complexe de groupe de Galois $G = \Gal(K/\mathbb Q)$
est une sous-extension de degré relatif au plus 2 du compositum 
d'une extension quadratique imaginaire $k = \mathbb Q[\sqrt{-d}]$ 
et d'une extension galoisienne (totalement) réelle de même groupe 
$G$.\medskip

\PreuveTh Notons $K_+$ le plus grand sous-corps (totalement) réel 
de $K$ et $C$ le groupe de Galois $\Gal(K/K_+)$. D'après ce qui 
précède, le groupe $C$ est d'ordre 2, engendré par la conjugaison 
complexe $\tau$ et distingué dans $G$, de sorte que $\tau$ est 
centrale dans $G$. Cela étant, pour toute extension quadratique 
totalement imaginaire $k = F[\sqrt{-d}]$ de $F$ qui n'est pas 
contenue dans $K$, l'extension composée $K' = Kk = K[\sqrt{-d}]$ 
est encore un corps à conjugaison complexe (comme composé de 
tels corps) et c'est aussi le produit direct sur $F$ de $K$ et de $k$, 
i.e. une extension galoisienne de groupe $\Gal(K'/F) \simeq G \times
C_2$. Mais c'est tout aussi bien le produit direct sur $F$ de $K'_+$ 
et de $k$, si $K'_+$ désigne le sous-corps réel maximal de $K$. Et 
on a donc $\Gal(K'_+/F) \simeq G$, comme annoncé.\medskip

\medskip
\noindent{\bf 1.b. Nombres de Weil usuels (i.e. à l'infini) et conjugaison complexe}\medskip

Notons toujours $\ov\QQ \subset \CC$ la clôture algébrique de $\QQ$ dans $\CC$ et écrivons $\tau$ la conjugaison complexe sur $\ov\QQ$, i.e. la restriction à $\ov\QQ$ de la conjugaison complexe.

\begin{Def}\label {W0}
On dit qu'un élément $a \in \ov\QQ$ est un nombre de Weil (à l'infini) lorsqu'on a~: $\mid a^\sigma \mid =1$ pour tout $\sigma$ de $G=\Gal(\ov\QQ/\QQ)$.\par
L'ensemble $W_{\ov\QQ}^\infty$ des nombres de Weil est ainsi le sous-groupe de $\ov\QQ^\times$~:\smallskip

\centerline{$W_{\ov\QQ}^\infty \ = \ \{a \in \ov\QQ^\times \mid \forall \sigma \in G  \quad a^{\sigma(1+\tau)}=1 \}$.}
\end{Def}

\Remarque Pour tout entier naturel $m$, on définit usuellement les {\it $m$-nombres de Weil} par la condition~: $a^{\sigma(1+\tau)}=m$, pour tout $\sigma$ de $G$~; en d'autres termes, ce sont les produits $\sqrt m \ b$, où $b$ est un nombre de Weil au sens de la définition~\ref{W0}.\medskip

Le lien entre nombres de Weil et conjugaison complexe se présente alors ainsi~:

\begin{Prop}
Soit $a \in W_{\ov\QQ}^\infty \setminus \{\pm 1\}$. Alors $K = \QQ[a]$ est un corps à conjugaison complexe~; autrement dit, il s'écrit $K = K_+[\sqrt{-\delta}]$, où $F$ est un corps
totalement réel et $\delta \gg 0$ un élément de $F$ totalement positif.
\end{Prop}

\Preuve De l'identité $a^{1+\tau}= 1$, on tire $a^{(1+\tau)\sigma} = 1 = a^{\sigma(1+\tau)}$~; d'où $a^{\sigma\tau}=a^{\tau\sigma}$ pour tout $\sigma$ de $G$. Soit alors $b=a+a^\tau \in \RR$. De $b^\sigma = a^\sigma + a^{\tau\sigma}= a + a^{\sigma\tau} \in \RR$, il suit que le sous-corps $F = \QQ[b]$ est totalement réel~; et de $a^\tau =a^{-1}$ que $a$ est racine du polynôme $P(X) = X^2 -bX +1$ de $F[X]$. Comme $a$ est non réel (en vertu de l'hypothèse $a \ne \pm 1$), il vient $\deg_Fa =2$ et $[\QQ[a] : \QQ[b]] = 2$.\par
Enfin, aucun des conjugués de $a$ n'étant réel, le corps $K=\QQ[a]$ est donc totalement imaginaire et le discriminant $-\delta$ du polynôme $P(X)$ est totalement négatif.

\begin{Cor}
Sous les mêmes hypothèses, la clôture galoisienne $\wh K$ de $K$ est encore un corps à conjugaison complexe et la restriction de la conjugaison $\tau$ à $\wh K$ est centrale dans $\Gal (\wh K/K)$. En particulier, $\wh K$ est une extension quadratique totalement imaginaire du corps galoisien (totalement) réel $\wh K_{tr} = \wh K^{<\tau>}$.
\end{Cor}

\Preuve Notons $\wh F$ la clôture galoisienne de $F$, i.e. le compositum des $F^\sigma$ lorsque $\sigma$ décrit $G=\Gal(\ov\QQ/\QQ)$. Nous avons $K=F[\sqrt{-\delta}]$ pour un $\delta \gg 0$ de $F$, donc $\wh K =\wh K_+[\sqrt{-\delta}]$, avec $\wh K_+=\wh F[\sqrt{\delta^\sigma/\delta},\ \sigma \in G]$, comme annoncé.

\begin{Th}
Notons $\ov\QQ_{tr}$ la plus grande sous-extension totalement réelle de $\ov\QQ$. On a alors~:

\centerline{$\QQ[W^\infty_{\ov\QQ}] \ = \ \ov\QQ_{tr}[i]$~;}
\noindent d'où~:\smallskip

\centerline{$W^\infty_{\ov\QQ} \ = \ \ov\QQ_{tr}[i]^{\times(\tau-1)} \ = \ \left\{ \frac{a-ib}{a+ib} \mid a\in \ov\QQ_{tr},\ b\in \ov\QQ_{tr},\ (a,b) \ne (0,0) \right\}$.}
\end{Th}

\Preuve Cela résulte du théorème 90 de Hilbert appliqué à l'extension $\ov\QQ_{tr}[i]/\ov\QQ_{tr}$.

\begin{Cor}
Soit $K \subset \CC$ un corps de nombres, $K_{tr}$ sa plus grande sous-extension totalement réelle et $K_{cm} = K_{tr}[W^\infty_K]$ le sous-corps à conjugaison complexe engendré sur $K_{tr}$ par les nombres de Weil contenus dans $K$. On a alors l'alternative suivante~:
\begin{itemize}
\item[(i)] ou bien on a : $[K_{cm} : K_{tr}] = 1$, auquel cas le groupe $W_K^\infty$ se réduit à $\{\pm 1\}$~;
\item[(ii)] ou bien on a : $[K_{cm} : K_{tr}] = 2$, auquel cas le corps $K_{cm}$ est à conjugaison complexe et il vient~:

\centerline{$W_K^\infty = W_{K_{cm}}^\infty = K_{cm}^{\times (\tau-1)}$.}
\end{itemize}
\end{Cor}

\Preuve Puisque $K$ contient $K_{cm}$, son sous-corps $K_{tr}$ est encore la plus grande sous-extension totalement réelle de $K_{cm}$. Cela étant, dans le premier cas, il suit directement~: $W_K^\infty=W_{K_{cm}}^\infty=W_{K_{tr}}^\infty = \{\pm 1\}$~; et dans le second, il vient, en vertu du théorème 90~:
$W_K^\infty = W_{K_{cm}}^\infty = K_{cm}^{\times (\tau-1)}$.

\bigskip\smallskip
\noindent{\large\bf 2. Corps à conjugaison $\ell$-adique et $\ell$-nombres de Weil}
\medskip

Nous allons maintenant définir l'analogue $\ell$-adique de la conjugaison complexe et des nombres de Weil pour chaque nombre premier donné $\ell$.

\medskip
\noindent{\bf 2.a. Corps totalement $\ell$-adiques et corps à conjugaison $\ell$-adique}\medskip

Nous supposons fixés dans cette section un plongement de la clôture algébrique $\ov\QQ$ de $\QQ$ dans le complété $\CC_\ell$ de la clôture algébrique $\ov\QQ_\ell$ du corps $\Ql$~; en d'autres termes, nous regardons $\ov\QQ$ comme un sous-corps de $\ov\QQ_\ell$. Cel point est essentiel pour parler sans ambiguïté de conjugaison $\ell$-adique.

\begin{Def}
Nous disons qu'un corps de nombres $K \subset \ov\QQ_\ell$ est~:
\begin{itemize}
\item[(i)] {\em totalement $\ell$-adique}, lorsque tous ses conjugués $K^\sigma$ par $\Gal(\ov\QQ/\QQ)$  sont contenus dans $\Ql$~; en d'autres termes lorsque toutes les places $\l$ de $K$ au-dessus de $\ell$ sont de degré $[K_\l : \Ql ]=1$~;
\item[(ii)] {\em à conjugaison $\ell$-adique}, lorsque c'est une extension non décomposée (aux places) au-dessus de $\ell$ d'un sous-corps $F$ totalement $\ell$-adique~;
\item[(iii)] {\em à conjugaison $\ell$-adique quadratique non ramifiée}, lorsque c'est une extension quadratique totalement inerte (aux places) au-dessus de $\ell$ d'un sous-corps $F$ totalement $\ell$-adique.
\end{itemize}
\end{Def}

Précisons tout cela à la lumière de la situation à l'infini~:
\begin{itemize}
\item[(i)] La première condition affirme qu'un corps $F$ est {\em totalement $\ell$-adique}, dès lors qu'on a l'égalité entre complétés en chaque place $\l$ de $F$ au-dessus de $\ell$~: \smallskip

\centerline{$F_\l =\Ql $.}\smallskip

C'est l'analogue exact de la condition qui définit un corps {\em totalement réel}.

\item[(ii)] La seconde condition affirme qu'un corps $K$ est {\em à conjugaison $\ell$-adique} dès lorsqu'il contient un sous-corps $F$ {\em totalement $\ell$-adique} tel qu'en chaque place $\l$ de $F$ au-dessus de $\ell$ le corps $K$ admette une unique place $\L$ au-dessus de $\l$~; de sorte que l'on a dans ce cas l'égalité entre degré global et degré local~:\smallskip

\centerline{$[K : F]  \  = \  [K_\L : F_\l] \  = \ [K_\L :\Ql]$.}

C'est évidemment un analogue  de la condition de {\em conjugaison complexe}.

\item[(iii)] La troisième condition dit que $K$ est {\em à conjugaison $\ell$-adique quadratique non ramifiée} lorsque les extensions locales $K_\L/F_\l$ sont non ramifiées et de degré 2. Elle demande donc une analogie plus étroite avec la situation classique  où l'extension locale à l'infini $\CC/\RR$ est naturellement de degré 2 et non ramifiée\footnote{En toute rigueur (cf. \cite{J3}), la complexification d'une place réelle n'est pas de la ramification.}. On a alors $F_\l =\Ql$ (puisque F est totalement $\ell$-adique) et $K_\L = \Ql [\zeta]$, où $\zeta$ désigne une racine primitive $(\ell+1)$-ème de l'unité dans $\ov\Ql$ (puisque $\Ql$ possède une unique extension quadratique non ramifiée, à savoir $\Ql [\zeta]$).
\end{itemize}
\smallskip

La condition de de degré et de non ramification introduite dans la définition de la conjugaison $\ell$-adique n'est pas anodine, comme le montre la proposition suivante~:

\begin{Prop} \label{CE}
On a les règles de composition suivantes~:
\begin{itemize}
\item[(i)] Le compositum de deux corps totalement $\ell$-adiques est totalement $\ell$-adique.
\item[(ii)] En revanche, le compositum de deux corps à conjugaison $\ell$-adique n'est pas nécessairement un corps à conjugaison $\ell$-adique.
\item[(iii)] Mais le compositum de deux corps à conjugaison $\ell$-adique quadratique non ramifiée est encore un corps à conjugaison $\ell$-adique quadratique non ramifiée.
\end{itemize}
\end{Prop}

\Preuve Examinons successivement les diverses assertions~:
\begin{itemize}
\item[(i)]
Si $F$ et $F'$ sont deux corps totalement $\ell$-adiques, on a par définition~: $F^\sigma \subset \Ql$ et $F'^\sigma \subset \Ql$ pour chaque $\sigma$ de $\Gal(\ov\QQ/\QQ)$~; donc $(FF')^\sigma = F^\sigma F'^\sigma\subset \Ql$, de sorte que le compositum $FF'$ est encore totalement $\ell$-adique.

\item[(ii)] Considérons, par exemple, deux extension absolument galoisiennes $K$ et $K'$, de groupes respectifs $G=\Gal(K/\QQ)$ et  $G'=\Gal(K'/\QQ)$ , linéairement disjointes sur $\QQ$ et non décomposées aux places au dessus de $\ell$. Sous ces hypothèses, le compositum $N=KK'$ est encore galoisien et son groupe de Galois est le produit direct $\Gal(N/\QQ) = G \times G'$. Les extensions locales $K_\L/\Ql$ et $K_\L/\Ql$  sont aussi galoisiennes de groupes respectifs $G=\Gal(K_\L/\Ql)$ et  $G'=\Gal(K'_\L/\Ql)$~; mais elles ne sont plus en général disjointes, de sorte que le sous-groupe de décomposition de la place $\L$ dans l'extension $N/\QQ$, qui s'identifie au groupe de Galois $\Gal(N_\L/\Ql)$, est alors le produit amalgamé~:\smallskip

\centerline{$D = \{(\sigma,\sigma')\in G\times G' \mid \sigma |_{K_\L \cap K'_\L} = \sigma'|_{K_\L \cap K'_\L}\}$.}

\noindent Il suit de là que le sous-groupe $D$ n'est normal dans $G \times G'$ que si l'extension locale $(K_\L \cap K'_\L)/\Ql$ est abélienne, ce qui n'est pas le cas en général~: $K$ et $K'$ sont alors à conjugaison $\ell$-adique, mais non leur compositum $KK'$.

\item[(iii)] Soient maintenant $K$ et $K'$ deux corps à conjugaison $\ell$-adique, tous deux extensions quadratiques totalement inertes (aux places au-dessus de $\ell$) de leurs sous-corps totalement $\ell$-adiques maximaux respectifs $F$ et $F'$. Dans le diagramme d'extensions~:

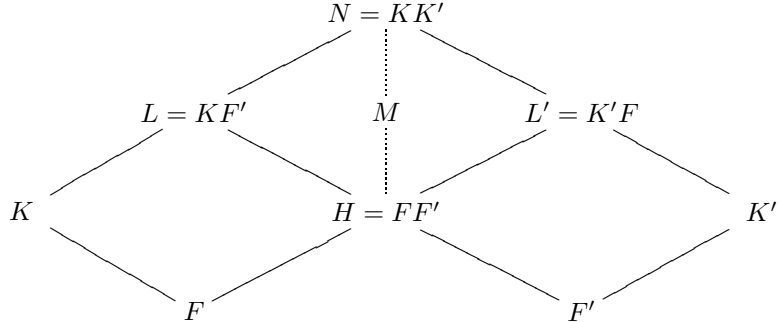
\begin{displaymath}
\xymatrix{
{} & {} & N=KK'\ar@{-}[ld] \ar@{.}[d] \ar@{-}[rd] & {} & {}  \\
{} & L=KF' \ar@{-}[rd] & M \ar@{.}[d] & L'=K'F \ar@{-}[ld]  & {} \\
\quad K \quad \ar@{-}[rd]  \ar@{-}[ru] &{} & H=FF' & {} & \quad K' \quad \ar@{-}[ld]  \ar@{-}[lu] \\
{} & F {} \ar@{-}[ru] & {} & F'  \ar@{-}[lu] & {}
}
\end{displaymath}\smallskip

le compositum $FF'$ est ainsi totalement $\ell$-adique et les extensions quadratiques $KF'/FF'$ et $K'F/FF'$  sont totalement inertes (aux places $\ell$-adiques).\smallskip

- Si elles coïncident, $N=KK'$ est une extension quadratique totalement inerte (aux places au-dessus de $\ell$) du corps totalement $\ell$-adique $H=FF'$.\par

-  Sinon, $N$ est biquadratique et non ramifiée sur $H$ (aux places au-dessus de $\ell$), mais non totalement inerte (puisque non cyclique), de sorte  que l'unique sous-extension quadratique $M$ de $N$ contenant $H$ autre que $L$ et $L'$ se trouve être le sous-corps de décomposition commun à toutes les places $\ell$-adiques~: en d'autres termes, $M$ est totalement $\ell$-adique et $N$ en est une extension quadratique totalement inerte (aux places au-dessus de $\ell$).
\end{itemize}

\medskip
\noindent {\bf 2.b. Définition des $\ell$-nombres de Weil (i.e. au sens $\ell$-adique)}
\medskip

D'après la définition \ref{W0}, les nombres de Weil classiques sont caractérisés à l'aide de la norme locale à l'infini $\nu_\infty = 1 + \tau = N_{\CC/\RR}$ par l'identité~:\smallskip

\centerline{$W_{\ov\QQ}^\infty \ = \ \{a \in \ov\QQ^\times \mid \forall \sigma \in G = \Gal(\ov\QQ/\QQ) \quad a^{\sigma(1+\tau)}=1 \}$.}\smallskip

Une autre façon d'exprimer cette même identité consiste à introduire le corps des racines $K_a$ du polynôme minimal $\Pi_a$ de $a$ sur $\QQ$, puis son sous-corps réel $H_a =K_a \cap \RR$, et à définir la norme locale $\nu_a$ comme l'opérateur $N_{K_a/H_a}$,  de sorte qu'avec ces notations il vient $\nu_a = \nu_\infty$ ou $\nu_a = 1$, suivant que $K_a$ est réel ou non~; et, dans les deux cas~:\smallskip

\centerline{$W_{\ov\QQ}^\infty \ = \ \{a \in \ov\QQ^\times \mid \forall \sigma \in G = \Gal(\ov\QQ/\QQ) \quad a^{\sigma \nu_a} \in \mu_\RR \}$,}\smallskip

\noindent où $\mu_\RR =\{\pm 1\}$ est le groupe des racines de l'unité dans $\RR$.\medskip

Sous cette dernière forme, la traduction $\ell$-adique de la définition des nombres de Weil est immédiate : désignons toujours par $K_a$ le corps des racines du polynôme minimal de $a$ sur $\QQ$~; mais écrivons maintenant $H_a = K_a \cap\, \Ql$ son sous-corps $\ell$-adique et $\nu_a$ la norme locale correspondante $N_{K_a/H_a}$. Cela étant~:

\begin{Def}\label {Wl}
Nous disons qu'un élément $a \in \ov\QQ$ est un {\em nombre de Weil $\ell$-adique} (ou encore un  $\ell$-{\it nombre de Weil}) lorsque la norme locale de chacun de ses conjugués est une racine de l'unité, ce qui s'écrit~:\smallskip

\centerline{$W_{\ov\QQ}^{(\ell )}\ = \ \{a \in \ov\QQ^\times \mid  \forall \sigma \in G=\Gal(\ov\QQ/\QQ)  \quad a^{\sigma \nu_a} \in \mu_{\Ql} \}$,}\smallskip

\noindent où $\mu_{\Ql}$ désigne le groupe (d'ordre 2 ou $\ell-1$) des racines de l'unité de $\Ql$.\par
L'ensemble $W_{\ov\QQ}^{(\ell)}$ des $\ell$-nombres de Weil est encore un sous-groupe de $\ov\QQ^\times$.
\end{Def}

Si un $\ell$-nombre de Weil $a$ possède un plongement $\ell$-adique (i.e. si l'on a $a^\sigma \in \Ql$ pour un $\sigma \in G$), il résulte clairement de la définition que $a^\sigma$, et donc $a$ lui-même, est une racine de l'unité. En particulier~:

\begin{Prop}
Si $K \subset \ov\QQ_\ell$ est un corps de nombres qui possède au moins un plongement $\ell$-adique, autrement dit dont l'un au moins des conjugués $K^\sigma$ est contenu dans $\Ql$, le groupe $W_K^{(\ell)}$ des $\ell$-nombres de Weil contenus dans $K^\times$ se réduit au sous-groupe $\mu_K$ des racines de l'unité dans $K$.
\end{Prop}

Si maintenant $a \in W_{\ov\QQ}^{(\ell )} \setminus \mu_{\ov\QQ}$ est un $\ell$-nombre de Weil qui n'est pas une racine de l'unité, le contre-exemple donné dans la proposition \ref{CE} montre qu'il n'y a pas lieu d'attendre une caractérisation aussi simple que dans le cas classique, sauf à imposer de sévères restrictions. Par exemple~:

\begin{DProp}
Convenons de dire qu'un $\ell$-nombre de Weil $a$ est $\ell$-quadratique lorsque le corps $\QQ[a]$ est à conjugaison $\ell$-adique quadratique non ramifiée. Alors~:\begin{itemize}
\item[(i)] La clôture galoisienne $K_a$ du corps $\QQ[a]$ engendré par un nombre de Weil $\ell$-quadratique est un corps à conjugaison $\ell$-adique quadratique non ramifiée.
\item[(ii)] Le sous-corps $\QQ_{\ell q}$ de $\ov\QQ$ engendré sur $\QQ$ par l'ensemble des nombres de Weil $\,\ell$-quadratiques est l'unique extension quadratique non ramifiée (aux places au-dessus de $\ell$) de l'extension totalement $\ell$-adique maximale $\QQ_{t\ell}$ de $\QQ$.
\end{itemize}
\end{DProp}

\Preuve Elle est identique {\it mutatis mutandis} à celle du cas classique. Si $a$ est un nombre de Weil $\ell$-adique qui est $\ell$ quadratique, la clôture galoisienne $K_a$ est bien à conjugaison $\ell$-adique quadratique , en vertu de la proposition \ref{CE} {\it (iii)}, comme composé de tels corps. Plus généralement, le compositum de tous les $K_a$ vérifie encore la même propriété. Et il reste simplement à voir que son sous-corps $\ell$-adique maximal est aussi l'extension totalement $\ell$-adique maximale $\QQ_{t\ell}$ de $\QQ$. Or~: \begin{itemize}
\item[(a)] Pour $\ell \ne 2$, notons $\zeta$ une racine primitive $2(\ell-1)$-ieme de l'unité dans $\ov\QQ$, et observons que c'est un nombre de Weil $\ell$-quadratique, de sorte que $\QQ_{\ell q}$ contient le corps cyclotomique $\QQ[\zeta]$. Il suit de là que, pour tout nombre algébrique $a$ totalement $\ell$-adique, le compositum $\QQ[\zeta^2,a]$ est encore toatalement $\ell$ -adique. Ainsi, le produit $\zeta a$ étant bien $\ell$-quadratique, $\QQ_{\ell q}$ contient $\QQ[a]$.
\item[(b)] Et pour $\ell = 2$, on conclut de même en prenant une racine cubique de l'unité.
\end{itemize}

\bigskip
\noindent {\large\bf 3. Nombres de Weil et valeurs absolues $\ell$-adiques}

\medskip

Commençons par rappeler ou préciser quelques résultats sur les valeurs absolues.

\medskip
\noindent {\bf 3.a. Comparaison des valeurs absolues réelles et $\ell$-adiques}
\medskip

Classiquement les valeurs absolues normalisées  réelles (i.e. à valeurs dans $\mathbb R_+$)  sont définies sur le groupe multilicatif d'un corps de nombres $K$ par les formules~:
$$
\mid x\mid^\infty _{\p} \ = \left \{ \aligned
& N{\p}^{-v_{\p} (x)}\  \ \ \ \ \  \ \ \ \  \ &\text{pour}\
{\p}\nmid \infty\, , \\
& \mid N_{K_{\p} /\RR} (x) \mid\ &\text{pour}\ {\p} \mid \infty\,;
\endaligned
\right.
$$
et la normalisation est dictée par la formule du produit~: $\prod_{\p} \mid x \mid^\infty_\p =1$.\smallskip

Par composition avec le logarithme réel, on obtient ainsi le plongement logarithmique usuel (où $\p$ décrit l'ensemble $Pl_K$ des places finies ou infinies de $K$)~:
\begin{displaymath}
\xymatrix{K^\times \ar[r]^{(\mid\ \mid^\infty_\p)_\p\ \quad \ }  & \oplus_{\p \in Pl_K} \RR_+^\times \ar[r]^{\Log_\infty} &  \oplus_{\p \in Pl_K} \RR}
\end{displaymath}

\noindent Le noyau des valeurs absolues réelles attachées aux places finies~:\smallskip

\centerline{$E_K^\infty = \{x \in K^\times \mid \ \forall \p \nmid \infty \quad \mid x \mid_\p^\infty \ = 1\}$}\smallskip

\noindent est ainsi le groupe des unités usuelles (à l'infini) du corps $K$~;  celui des valeurs absolues attachées aux places à l'infini~:\smallskip

\centerline{$W_K^\infty = \{x \in K^\times \mid \ \forall \p \mid \infty \quad \mid x \mid_\p^\infty\ = 1\}$}\smallskip

\noindent est le groupe de Weil (à l'infini) défini plus haut~; enfin l'intersection~:\smallskip

\centerline{$\mu_K^\infty \ = \ E_K^\infty \cap W_K^\infty $,}\smallskip

\noindent i.e. le noyau de toutes les valeurs absolues réelles, se réduit au sous-groupe $\mu_K^\infty$ de $K^\times$ formé des racines de l'unité du corps $K$, d'après un vieux résultat de Kronecker.\medskip

En ce qui concerne les valeurs absolues $\ell$-adiques (i.e. à valeurs $\ell$-adiques), les choses sont plus compliquées~: classiquement les valeurs absolues {\it principales} normalisées sont définies sur le groupe multilicatif d'un corps de nombres $K$ par les formules~:
$$
\mid x\mid^{\langle\ell\rangle} _{\p} \ = \left \{ \aligned
& \ 1 &\text{pour}  & \ {\p}\mid \infty\, , \\ 
& \langle N{\p}^{-v_{\p} (x)} \rangle\  \ \ \ \ \  \ \ \ \  \ &\text{pour} & \ {\p}\nmid \ell\infty\, , \\
& \langle N_{K_{\p} /\Ql} (x) N{\p}^{-v_{\p} (x)}\rangle\ &\text{pour} & \ {\p} \mid \ell\ ;
\endaligned
\right.
$$
où $\langle \cdot \rangle$ désigne la projection naturelle du groupe multiplicatif $\Zl^\times$ sur son sous-groupe sans torsion $\langle\Zl^\times\rangle=1+2\ell\Zl$~; et la normalisation est dictée par la formule du produit~: $\prod_{\p} \mid x \mid^{\langle\ell\rangle}_\p =1$. C'est le point de vue adopté dans \cite{J2,J3}.\smallskip

Par composition avec le logarithme $\ell$-adique, on obtient ainsi le plongement logarithmique (où $\p$ décrit maintenant l'ensemble $Pl_K^{\,0}$ des places finies de $K$)~:
\begin{displaymath}
\xymatrix{K^\times \ar[r]^{(\mid\ \mid^{\langle\ell\rangle} _\p)_\p \quad \quad }  & \oplus_{\p \in Pl^{\,0}_K} \langle\Zl^\times\rangle  \ar[r]^{\quad \Log_\ell\quad } &  \oplus_{\p \in Pl^{\,0}_K}\, \Ql}
\end{displaymath}

Malheureusement, cette définition n'est pas suffisante pour rendre compte des propriétés normiques des $\ell$-nombres de Weil ; et c'est pourquoi, il est préférable de considérer dans ce contexte les valeurs absolues {\it non principales} à valeurs dans le groupe multiplicatif sans torsion $\langle\Zl^\times\rangle \; \ell^{\,\QQ}$, en procédant comme suit~: faisons choix d'un plongement $\phi_\ell$ du corps des nombres complexes $\CC$ dans la clôture algébrique $\CC_\ell$ du complété $\ell$-adique $\QQ_\ell$ de $\QQ$. Cela étant, à chaque plongement à l'infini de $K^\times$\smallskip

\centerline{$x \mapsto x_\p$}\smallskip

\noindent dans $\RR^\times$ ou $\CC^\times$, correspond alors un plongement de $K^\times$ dans $\CC_\ell^\times$, donné par\smallskip

\centerline{$x \mapsto \phi_\ell(x_\p)$.}

\begin{Def}
La lettre $\p$ désignant maintenant  soit une place finie soit un plongement à l'infini d'un corps de nombres $K$, nous appelons {\em valeur absolue $\ell$-adique (non principale)} attachée à $\p$ l'application à valeurs dans le groupe sans torsion $(1+2\ell\Zl)\,\ell^{\,\QQ}$ qui est définie sur le groupe multiplicatif $K^\times$ par la formule~:
$$
\mid x\mid^{(\ell)} _{\p} \ = \left \{ \aligned
& \ \ell^{-v_\ell(\phi_\ell(x_\p))} &\text{pour}  & \ {\p}\mid \infty\, , \\ 
& \langle N{\p}^{-v_{\p} (x)} \rangle\  \ \ \ \ \  \ \ \ \  \ &\text{pour} & \ {\p}\nmid \ell\infty\, , \\
& \langle N_{K_{\p} /\Ql} (x) N{\p}^{-v_{\p} (x)}\rangle \ N{\p}^{v_{\p} (x)} &\text{pour} & \ {\p} \mid \ell\ ;
\endaligned
\right.
$$
\end{Def}

\Remarque La définition de la {\it valeur absolue $\ell$-adique} $| \,x\,|^{(\ell)} _{\p}$ en une place à l'infini  $\p$ dépend évidemment du choix du plongement $\phi_\ell$ de $\CC$ dans $\CC_\ell$~; en revanche, la famille $(\mid x\mid^{(\ell)} _{\p})_{\p\mid\infty}$, est, à l'ordre près, indépendante du choix effectué.

\begin{Prop}
Les valeurs absolues non principales vérifient encore la formule du produit~:

\centerline{$\prod_{\p} \mid x \mid^{(\ell)}_\p =1$.}
\end{Prop}

\Preuve Un calcul élémentaire donne, en effet~:\smallskip

\centerline{$\prod_{\p \mid \infty} \mid x \mid^{(\ell)}_\p \;=\; \ell^{-v_\ell(\phi_\ell(N_{K/\QQ}(x)}\;=\; \ell^{-v_\ell(N_{K/\QQ}(x)}\;=\; \prod_{\p \mid \ell}N\p^{-v_\p(x)}$ ;}\smallskip

\noindent d'où le résultat annoncé, le produit des valeurs absolues principales étant égal à 1.

\begin{Prop} \label{mu}
Avec les conventions précédentes, il vient :
\begin{itemize}
\item[(i)] Le noyau des valeurs absolues $\ell$-adiques attachées aux  places $\p$ qui ne divisent pas $\ell$ est encore le groupe des unités usuelles du corps $K$~:\smallskip

\centerline{$E_K^{(\ell)} = \{x \in K^\times \mid \ \forall \p \nmid \ell \quad \mid x \mid_\p^{(\ell)} \ = 1\} = E_K^\infty$.}

\item[(ii)] Le noyau des valeurs absolues $\ell$-adiques attachées aux places  $\p$ qui divisent $\ell$ est le $\ell$-groupe de Weil défini plus haut~:\smallskip

\centerline{$W_K^{(\ell)} = \{x \in K^\times \mid \ \forall \p \mid \infty \quad \mid x \mid_\p^{(\ell)}\ = 1\}$.}

\item[(iii)] Le noyau de toutes les valeurs absolues $\ell$-adiques, i.e. l'intersection\smallskip

\centerline{$\mu^{(\ell)}_K = E_K^{(\ell)} \cap W_K^{(\ell)}\ =  \{x \in K^\times \mid \ \forall \p \in Pl_K  \quad \mid x \mid_\p^{(\ell)} \ = 1\}$.}\smallskip

\noindent est le sous-groupe de $K^\times$ formé des unités (au sens usuel) qui sont normes à chaque étage fini de la $\Zl$-extension cyclotomique $K^{(\ell)}_\infty$ du corps $K$.
\end{itemize}
\end{Prop}

\Preuve La condition ($i$) impose $v_\p(x) =1$ à la fois pour $\p \mid \ell$ et pour $\p \nmid \ell\infty$~;  de sorte que les éléments de $E^{(\ell)}_K $ sont encore les unités (au sens ordinaire) du corps $K$. \smallskip

En termes normiques, la condition ($ii$) s'écrit~: $N_{K_{\p} /\Ql} (x) \in \mu_{\Ql}, \ \forall \p \mid \ell$, où $\mu_{\Ql}$ désigne le sous-groupe de torsion de $\Ql$~; elle caractérise donc les $\ell$-nombres de Weil.\smallskip

Enfin, la condition normique $\mid x\mid^{\langle\ell\rangle} _{\p} = \langle N_{K_{\p} /\Ql} (x) N{\p}^{-v_{\p} (x)}\rangle = 1,\ \forall \p \mid \ell$ caractérise précisément les éléments qui sont normes locales aux places au-dessus de $\ell$ dans la $\Zl$-extension cyclotomique $K^{(\ell)}_\infty$ de $K$ (cf. e.g. \cite{J2,J3}). D'où le résultat annoncé, en vertu du principe de Hasse, puisque les unités sont évidemment normes locales aux places étrangères à $\ell$ dans l'extension procyclique $K^{(\ell)}_\infty/K$.\smallskip

Nous allons voir que ce groupe $\mu^{(\ell)}_K$ ne se réduit pas toujours au sous-groupe $\mu^\infty_K$.

\medskip
\noindent {\bf 3.b. \'Etude du noyau des valeurs absolues $\ell$-adiques}
\medskip

D'après ce qui précède, le goupe $\mu_K^{(\ell)}$ est un sous-groupe du groupe des unités $E_K^\infty$ qui contient le sous-groupe des racines de l'unité $\mu_K^\infty$. Dans la double inclusion \smallskip

\centerline{$ \mu_K^\infty \ \subset \ \mu_K^{(\ell)} \ \subset \ E^\infty_K$}\smallskip

\noindent il est facile de voir que l'égalité est, en fait, possible à droite comme à gauche, comme l'illustrent les propositions suivantes~:

\begin{Prop}
Si $K$ est un corps de nombres totalement $\ell$-adique, le noyau $\mu^{(\ell)}_K$ des valeurs absolues $\ell$-adiques se réduit au groupe $\mu^\infty_K$ des racines de l'unité.
\end{Prop}

\Preuve C'est immédiat car, la norme locale $\nu_K = N_{K/(K\cap\Ql)}$  étant alors l'identité, on a immédiatement~:\smallskip

\centerline{$\mu_K^{(\ell)} \ = \ \{a \in E^\infty_K \mid a^{\nu_K} \in \mu_{\Ql} \}\ = \ \{a \in E^\infty_K \mid a \in \mu_{\Ql} \} \ = \ \mu^\infty_K$.}

\begin{Prop}
Si $K$ est un corps à conjugaison $\ell$-adique admettant $\QQ$ comme sous-corps totalement $\ell$-adique maximal (en d'autres termes, si $K$ est un corps de nombres dans lequel aucune des places $\ell$-adiques n'admet de décomposition), le noyau $\mu^{(\ell)}_K$ des valeurs absolues $\ell$-adiques dans  $K$ est égal au groupe des unités $E^\infty_K$.
\end{Prop}

\Preuve Dans ce cas, en effet, la norme locale $\nu_K = N_{K/(K\cap\Ql)}$ coïncide avec la norme globale $N_{K/\QQ}$~; et il vient donc~:

\centerline{$\mu_K^{(\ell)} \ = \ \{a \in E_K \mid a^{\nu_K} \in \mu_{\Ql} \}\ = \ \{a \in E_K \mid N_{K/\QQ}(a) \in \mu_{\QQ} \} \ = \ E_K$,}\smallskip

\noindent puisque les unités d'un corps de nombres sont évidemment de norme $\pm 1$.\medskip

Plus généralement, la détermination du groupe $\mu_K^{(\ell)}$ peut être menée à bien très simplement dès lors que le corps considéré est à conjugaison $\ell$-adique~: 

\begin{Th}
Soit K une extension galoisienne de degré $d$ non décomposée (aux places au-dessus de $\ell$) d'un sous-corps $F$ totalement $\ell$-adique~; pour chaque place $\p$ de $F$, désignons par $\chi_\p = \Ind_{D_\p}^G \ 1_{D_\p}$ l'induit à $G=\Gal(K/F)$ du caractère de la représentation unité du sous-groupe de décomposition de l'une quelconque des places de $K$ au-dessus de $\p$. \par
Avec ces notations, le caractère du $\QQ[G]$-module multiplicatif $\ \QQ \otimes_\ZZ \mu_K^{(\ell)}$ est donné par la formule~:
$$
\chi_{\mu^{(\ell)}_K} = \sum_{\p_{\!\infty} \in Pl^\infty_K}\, ( \chi_{\p_{\!\infty}} -1)
$$
\end{Th}

\Preuve On sait par le théorème de représentation de Herbrand que le tensorisé multiplicatif $\QQ \otimes_\ZZ E_K^\infty$ du groupe des unités est un $\QQ[G]$-module de caractère~: 

\centerline{$\chi_\infty =(\sum_{\p_{\!\infty} \in Pl^\infty_K}\  \chi_{\p_{\!\infty}}) - 1=\sum_{\p_{\!\infty} \in Pl^\infty_K}\  (\chi_{\p_{\!\infty}} - 1) + (\mid Pl_K^\infty \mid -1)\,1$.}

En d'autres termes, c'est la somme directe de $|\,Pl_K^\infty \,|$ sous-modules monogènes (notés multiplicativement) $M_{\p_{\!\infty}} = (1\otimes x_{\p_{\!\infty}})^{\QQ[G]}$ de caractères respectifs $(\chi_{\p_{\!\infty}} - 1)$ et de $(|\, Pl_K^\infty \,| -1)$ sous-modules de caractère unité. Ce point acquis, le noyau de la norme locale $\nu =\sum_{\sigma \in G} \sigma$, autrement dit le noyau de l'idempotent $\nu /d$, est l'image de l'idempotent supplémentaire $1-\nu/d$~; et c'est donc la somme directe des  $\mid Pl_K^\infty \mid$ modules $M_{\p_{\!\infty}}$. D'où la formule annoncée~: $\chi_\infty =\sum_{\p_{\!\infty} \in Pl^\infty_K} \,(\chi_{\p_{\!\infty}}- 1)$

\begin{Sco} \label{Sco}
Soit $K$ un corps de nombres  à conjugaison $\ell$-adique et $F$ son sous-corps totalement $\ell$-adique maximal.  Le noyau $\mu_K^{(\ell)}$ des valeurs absolues $\ell$-adiques dans $K^\times$ est alors la somme directe du sous-groupe fini $\mu_K^\infty$ et d'un $\ZZ$-module libre de dimension~:

\centerline{$\rg_\ZZ \mu^{(\ell)}_K = \rg_\ZZ E_K^\infty -\rg_\ZZ E_F^\infty = (r_K+c_K)-(r_F+c_F) $,}\smallskip

\noindent où $r$ et $c$ dénombrent les places réelles ou complexes des corps $K$ et $F$ considérés.
\end{Sco}

\Preuve D'après le théorème, lorsque $K$ est galoisien sur $F$ de groupe $G$, le $\ZZ[G]$-module $\mu_K^{(\ell)}$ est la somme directe de son sous-module de torsion $\mu_K^{\infty}$ et d'un $\ZZ$-module libre de dimension~: $\deg \chi_{\mu^{(\ell)}_K}=\sum_{\p_{\!\infty} \in Pl^\infty_K}\ (\deg \chi_{\p_{\!\infty}} -1)=(r_K+c_K)-(r_F+c_F)$.\par
Il s'agit de voir que la formule de rang obtenue est indépendante de l'hypothèse galoisienne utilisée pour la démontrer. Or, si $d$ désigne le degré de l'extension $K/F$ et $\nu$ la norme associée, le tensorisé $\QQ \otimes_\ZZ \mu_K^{(\ell)}$ n'est autre que le noyau de l'idempotent $\nu/d$ dans le $\QQ$-espace $\QQ \otimes_\ZZ E_K^\infty$ construit sur les unités de $K$~; et son image est le sous-espace $\QQ \otimes_\ZZ E_F^\infty$ construit sur les unités de $F$.

\begin{Cor}
Sous les mêmes hypothèses que dans le scolie~:
\begin{itemize}
\item[(i)] L'égalité $\mu_K^{(\ell)} = \mu_K^{\infty}$ a lieu si et seulement si le corps $K$ est soit totalement $\ell$-adique, soit extension quadratique totalement imaginaire et à conjugaison  $\ell$-adique d'un sous-corps $F$ à la fois totalement réel et totalement $\ell$-adique.
\item[(ii)] L'égalité $\mu_K^{(\ell)} = E_K^{\infty}$ a lieu si et seulement si le sous-corps totalement $\ell$-adique $F$ est le corps des rationnels $\QQ$ ou un corps quadratique imaginaire $\QQ[\sqrt{-\delta}]$.
\end{itemize}
\end{Cor} 

\Preuve La première égalité s'écrit~: $\rg_\ZZ E^\infty_K = \rg_\ZZ E^\infty_F$~; et la seconde~: $\rg_\ZZ E^\infty_F=0$.

\bigskip\medskip
\centerline{\bf \large Appendice : lien avec les unités logarithmiques}
\bigskip

D'après la proposition \ref{mu}, le noyau $\mu_K^{(\ell)}$ de toutes les valeurs absolues est le sous-groupe des unités du corps $K$ qui sont normes à chaque étage fini de la $\Zl$-extension cyclotomique $K_\infty^{(\ell)}/K$.\smallskip

Or, dans la théorie $\ell$-adique du corps de classes (cf. \cite{J2,J3}), le groupe des $\Zl$-normes cyclotomiques globales est le sous-groupe des unités logarithmiques $\wi\E_K$ du tensorisé multiplicatif $\,\R_K = \Zl \otimes_\ZZ K^\times$, qui est défini à partir des valuations logarithmiques par~:

\centerline{$\wi\E_K \ = \ \{x \in \R_K \mid \wi v_\p(x) = 0 \quad \forall \p \in Pl^{\,0}_K \}$.}\smallskip

Rassemblant alors la description de $\mu_K^{(\ell)}$ donnée par la proposition \ref{mu} et le calcul de rang effectué dans le scolie \ref{Sco}, nous obtenons immédiatement~:

\begin{Th}
Soient $K \subset \ov \QQ_\ell$ un corps de nombres à conjugaison $\ell$-adique, puis $F = K \cap \Ql$ son sous-corps totalement $\ell$-adique maximal, et $\nu = N_{K/F}$ la norme locale associée à $\ell$. Notons $\,\R_K = \Zl \otimes_\ZZ K^\times$ le $\ell$-groupe des idèles principaux du corps $K$, puis $\,\E_K = \Zl \otimes_\ZZ E_K^\infty$ le sous-groupe unité de $\R_K$, et $\wi\E_K$ le groupe des unités logarithmiques. Avec ces notations, le tensorisé $\ell$-adique $\V_K$ du noyau $\mu_K^{(\ell)}$ est l'intersection~:

\centerline{$\V_K := \Zl \otimes_\ZZ\mu_K^{(\ell)} = \{\varepsilon \in \E_K  \mid \varepsilon^\nu =1 \}\ = \ \E_K \cap\, \wi\E_K$.}\medskip

\noindent Et c'est le produit direct du $\ell$-groupe $\Zl \otimes_\ZZ \mu_K^\infty$ des racines de l'unité de $K$ par un $\Zl$-module libre de dimension~:\smallskip

\centerline{$\rg_{\Zl} \V_K =(r_K+c_K)-(r_F+c_F)$,}\smallskip

\noindent où $r$ et $c$ dénombrent les places réelles ou complexes des corps $K$ et $F$ considérés.
\end{Th}

\Remarque On prendra garde que le produit tensoriel par $\Zl$ a tué le groupe $\mu_{\Ql}$.

\begin{Cor}
Sous les hypothèses du théorème, le $\ell$-groupe $\wi \E_K$ des unités logarithmiques du corps $K$ est la racine dans $\R_K$ du produit $\wi \E_F \, \V_K$ du $\ell$-groupe $\wi \E_K$ des unités logarithmiques de $F$ par le noyau $\V_K$ des valeurs absolues (non principales)~:\smallskip

\centerline{$\wi\E_K =  \sqrt{\wi \E_F \, \V_K}$.}\smallskip

\noindent Il suit de là que le défaut $\delta_K$ de la conjecture de Gross généralisée pour le corps $K$ (et le premier $\ell$) est le même que celui $\delta_F$ de $F$. En particulier, $K$ satisfait la conjecture de Gross généralisée si et seulement si son sous-corps $F$ la satisfait aussi.
\end{Cor}

\Preuve Nous avons trivialement~: $\wi \E_K^\nu \subset \wi \E_F$~; donc, en notant $d = [K:F]$ le degré de l'extension $K/F$, la double inclusion~: \smallskip

\centerline{$\wi \E_K^d \subset \wi \E_F \, \V_K \subset \wi \E_K$~;}\smallskip

\noindent ce qui nous donne le résultat annoncé, les unités logarithmiques constituant un sous-module pur de $\R_K$. Et comme l'intersection $\wi\E_F \cap \V_F$ se réduit au sous-groupe $\Zl\otimes_\ZZ\mu^\infty_K$ des racines de l'unité dans $\R_K$, nous en tirons la formule de rang~:\smallskip

\centerline{$\rg_{\Zl} \wi\E_K = \rg_{\Zl} \wi\E_F + \rg_{\Zl} \V_K$.}\smallskip

\noindent Cela étant, si maintenant $\delta_K$ et $\delta_F$ désignent respectivement le défaut de la conjecture de Gross dans $K$ et dans $F$ (cf. e.g. \cite{FG,G,J1}), il vient ainsi~:\smallskip

\centerline{$r_K + c_K + \delta_K = (r_F + c_F + \delta_F) + [(r_K + c_K) - (r_F + c_F) ] = r_K + c_K + \delta_F$~;}\smallskip

\noindent d'où l'égalité annoncée~:

\centerline{$\delta_K = \delta_F$.}\medskip

\Remarque Le corps $K$ satisfait ainsi la conjecture de Gross généralisée (pour le nombre premier $\ell$) dès que son sous-corps $F$ est abélien sur $\QQ$, alors même qu'il peut très bien ne pas être lui-même absolument galoisien.

\Exemple Un corps de nombres $K$ absolument galoisien à groupe de Galois diédral vérifie ainsi la conjecture de Gross généralisée en tout premier $\ell$ qui ne se décompose pas au-dessus de son unique sous-corps quadratique $F$.

\begin{Th}
Soient $K \subset \ov \QQ_\ell$ un corps de nombres à conjugaison $\ell$-adique~; $F = K \cap \Ql$ son sous-corps totalement $\ell$-adique maximal~; et $c_F$ le nombre de places complexes de $F$. Si $K$ satisfait la conjecture de Leopoldt (pour le nombre premier $\ell$), le groupe de Galois $\Gal(K^{cr}\!/K)$ de la plus grande pro-$\ell$-extension abélienne $K^{cr}$ de $K$ qui est cyclotomiquement ramifiée est un $\Zl$-module de rang~:\smallskip

\centerline{$\rg_{\Zl}\Gal(K^{cr}\!/K) = c_F +1$.}
\end{Th}

Ce résultat, qui ne présuppose aucune hypothèse sur le comportement des places à l'infini, fournit en particulier une minoration de l'invariant lambda d'Iwasawa attaché aux $\ell$-groupes de classes d'idéaux $\Cl_{K_n}$ dans la $\Zl$-extension cyclotomique $K^c = \cup_{n\in\NN} K_n$ du corps $K$ considéré~:

\begin{Cor}
Sous les hypothèses et avec les notations du théorème, l'invariant {\it lambda} d'Iwasawa attaché aux $\ell$-groupes de $\ell$-classes d'idéaux dans la $\Zl$-extension cyclotomique $K^c/K$ vérifie la minoration~:\smallskip

\centerline{$\lambda_K \ge c_F$.}
\end{Cor}

\PreuveTh Rappelons qu'une pro-$\ell$-extension abélienne $L$ d'un corps de nombres $K$ est dite cyclotomiquement ramifiée (cf. \cite{JS}, Déf. 1) lorsque~:
\begin{itemize}
\item[(i)] toutes les places au-dessus de $\ell$ se ramifient infiniment dans $L/K$ et que
\item[(ii)] l'extension induite $LK^c/K^c$ sur la tour cyclotomique $K^c$ est non ramifiée.
\end{itemize}
D'après la théorie $\ell$-adique globale du corps de classes (cf. \cite{JS}, Th. 2), le sous-groupe de normes qui fixe $K^{cr}$ est le produit $\wh\U_K\R_K$ où $\wh\U_K = \U_K \cap \wi\U_K$ est le groupe des unités locales qui sont normes cyclotomiques et $\R_K$ désigne, comme d'habitude, le sous-groupe principal du $\ell$-adifié $\,\J_K$ du groupe des idèles de $K$. Cela étant, pour évaluer le $\Zl$-rang du groupe $\Gal(K^{cr}\!/K) \simeq \J_K/\wh\U_K\R_K$, il est commode de le remplacer par le sous-groupe d'indice fini 
$\Gal(K^{cr}\!/K^{nr}) \simeq \U_K\R_K/\wh\U_K\R_K$, où $K^{nr}$ désigne le $\ell$-corps des classes de Hilbert de $K$, i.e. la plus grande $\ell$-extension abélienne du corps $K$ qui est non ramifiée. Il vient ainsi~:\smallskip

\centerline{$\Gal(K^{cr}\!/K)  \sim \U_K\R_K/\wh\U_K\R_K \simeq \U_\ell/\wh\U_\ell \E^\infty_\ell)$,}\smallskip

\noindent où $\U_\ell = \prod_{\l|\ell}\U_\l$ est le groupe des unités semi-locales attaché aux places $\ell$-adiques, $\wh\U_\ell$ est le noyau des valeurs absolues $\ell$-adiques dans $\U_\ell$, et $\E_\ell$ est l'image dans $\U_\ell$ du sous-groupe unité $\E_K = \Zl \otimes E^\infty_K$ de $\R_K$. \par
Si maintenant $K$ satisfait la conjecture de Leopoldt (pour le premier $\ell$), le groupe $\E_K$ s'injecte dans $\U_\ell$ et l'intersection $\E_\ell \cap\, \wh\U_\ell$ est l'image bijective du sous-groupe $\V_K = \E_K^\infty \cap \wi\E_K$. Il vient donc~:\smallskip

\centerline{$\rg_{\Zl}\Gal(K^{cr}\!/K) = \rg_{\Zl}\Gal(K^{cr}\!/K^{nr}) = \rg_{\Zl}\U_\ell - \rg_{\Zl} \wh\U_\ell - \rg_{\Zl} \E_K + \rg_{\Zl} \V_K$.}\smallskip

\noindent Or, sous l'hypothèse de conjugaison $\ell$-adique, $\wh\U_\ell$ est le noyau dans $\U_\ell$ de l'opérateur norme $\nu = N_{K/F}$, dont l'image est d'indice fini dans le sous-groupe correspondant  relatif à $F$, de $\Zl$-rang égal au degré $[F:\QQ] = r_F + 2 c_F$. En fin de compte, il suit~:\smallskip

\centerline{$\rg_{\Zl}\Gal(K^{cr}\!/K) = (r_F+2c_F) -(r_K+c_K-1)+[(r_K+c_K)-(r_F+c_F)]$~;}\smallskip

\noindent ce qui donne précisément la formule annoncée~: $\rg_{\Zl}\Gal(K^{cr}\!/K) = c_F +1$.


\def\refname{\small{\sc  Références}}

\bigskip\noindent
{\small
\begin{tabular}{l}
{Jean-Fran\c cois {\sc Jaulent}}\\
Institut de Math{\'e}matiques de Bordeaux \\
Université {\sc Bordeaux 1} \\
351, cours de la lib{\'e}ration\\
F-33405 {\sc Talence} Cedex\\
courriel : Jean-Francois.Jaulent@math.u-bordeaux1.fr 
\end{tabular}
}

 \end{document}